\documentclass[12pt]{article}
\usepackage{amsmath,amsthm,amsfonts,amssymb,graphicx,booktabs}
\usepackage{flafter}
\allowdisplaybreaks[1]
\newcommand{\Z}{\mathbb Z}
\renewcommand{\pmod}[1]{\,(\textup{mod}\,#1)}
\numberwithin{equation}{section}
\theoremstyle{plain}
\newtheorem{theorem}{Theorem}[section]
\newtheorem{lemma}[theorem]{Lemma}
\newtheorem{corollary}[theorem]{Corollary}

\begin{document}
\begin{center}{\textbf{\large Divisibility Properties of Odd Part and Even Part Sums in Cubic Partitions}}\end{center}

	\begin{center}
		{\bf Gaurab Bardhan$^1$ and Nipen Saikia$^{2, \ast}$}\\
		$^1$Department of Mathematics, Tyagbir Hem Baruah College,\\ Jamugurihat, Sonitpur, Assam, India.\\
		E-mail: gaurabbardhan561@gmail.com
		\vskip2mm
		$^2$Department of Mathematics, Rajiv Gandhi
		University,\\ Rono Hills, Doimukh, Arunachal Pradesh, India.\\
		E-mail: nipennak@yahoo.com\\
		$^\ast$\textit{Corresponding author}.\end{center}
\vspace{2mm}
\begin{abstract}A cubic partition of a positive integer $n$ is a partition of $n$ in which even parts may appear in two distinct colors. The purpose of this paper is to study the arithmetic functions $COS(n)$ and $CES(n)$, which denote the sums of all odd and even parts, respectively, in all cubic partitions of $n$. We obtain their generating functions and divisor-sum formulas, and establish Ramanujan-type congruences modulo $3$, $4$, $5$, and $7$. To prove these congruences, we use Lambert series, classical theta-function identities, dissections, coefficient extraction, and weighted Euler identities.
\end{abstract}

\vspace{0.8em}

\section{Introduction}

Throughout this paper, for complex numbers $a$, $b$, and $q$ with $|q|<1$, we use the standard notation 
\begin{equation}\label{e1}
(a;q)_\infty=\prod_{j\ge0}(1-aq^j),\qquad \text{and}\qquad
(a,b;q)_\infty=(a;q)_\infty(b;q)_\infty.\end{equation}
 For convenience, we write $f_r=(q^r;q^r)_\infty$, where $r$ is a positive integer.

A non-increasing finite sequence of positive integers $\beta_1\geq \beta_2 \geq \cdots \geq \beta_s>0$ is said to be a partition of a positive integer $n$ if $n=\sum_{i=1}^{s}\beta_i$.
	The integers $\beta_i$ are called the parts of the partition.  
	If $p(n)$ denotes the number of partitions of $n$, then its generating function \cite{euler1748introductio} is given by
	$$\sum_{n=0}^{\infty}p(n)q^n=\dfrac{1}{(q;q)_\infty}=\dfrac{1}{f_1};\qquad p(0)=1.$$
	For example, $p(7)=15$, and the corresponding partitions are $7, 6 + 1, 5 + 2, 5 + 1 + 1, 4 + 3, 4 + 2 + 1, 4 + 1 + 1 + 1, 3 + 3 + 1, 3 + 2 + 2, 3 + 2 + 1 + 1, 3 + 1 + 1 + 1 + 1, 2 + 2 + 2 + 1, 2 + 2 + 1 + 1 + 1, 2 + 1 + 1 + 1 + 1 + 1$, and  $1 + 1 + 1 + 1 + 1 + 1 + 1$.
	Among the most celebrated arithmetic properties of $p(n)$ are the following
 congruences due to Ramanujan \cite{Ramanujan}:
	$$
	p(5n+4)\equiv0\pmod5,\quad p(7n+5)\equiv0\pmod7,\quad p(11n+6)\equiv0\pmod{11},
	$$
	valid for every nonnegative integer $n$. Ramanujan's congruences have motivated extensive investigations of divisibility properties of several restricted partition functions.
    
Chan \cite{Chan2010} studied the cubic partition function $a(n)$, in which even parts may appear in two distinct colors. The generating function of $a(n)$ is given by
\begin{equation}\label{e2}
A(q):=\sum_{n\ge0}a(n)q^n
=\dfrac{1}{f_1f_2}
=\prod_{m\ge1}\dfrac{1}{(1-q^{2m-1})(1-q^{2m})^2};\qquad a(0)=1.
\end{equation} 
For example, the cubic partitions of $3$ are $3$, $2_r+1$, $2_g+1$, and $1+1+1$, where the subscripts $r$ and $g$ distinguish the two colors of the even part. Clearly, $a(3)=4$.
Chan \cite[Theorem~1]{Chan2010} proved the identity
\begin{equation}\label{e3}
\sum_{n\ge0}a(3n+2)q^n
=3\dfrac{f_3^3f_6^3}{f_1^4f_2^4},
\end{equation}
which immediately gives the Ramanujan-type congruence $a(3n+2)\equiv0\pmod3$. Chen and Lin \cite[Theorems~1.3--1.4]{ChenLin2009} proved that, for every $n\ge0$,
$$
a(25n+22)\equiv0\pmod5
$$
and
\begin{align*}
a(49n+15)&\equiv0\pmod7,\\
a(49n+29)&\equiv0\pmod7,\\
a(49n+36)&\equiv0\pmod7,\\
a(49n+43)&\equiv0\pmod7.
\end{align*}
Xiong \cite[Theorem~1.3]{Xiong2011} extended the modulo $5$ result to an infinite family modulo powers of $5$.

Recently, Andrews and Ghosh Dastidar \cite{AndrewsDastidar2026}
defined $SOME(n)$ as the difference between the sum of all odd
parts and the sum of all even parts in the partitions of a
nonnegative integer $n$. They established several divisibility
properties of $SOME(n)$. In particular, they proved \cite[Theorem~4]{AndrewsDastidar2026} that
$SOME(5n+4)\equiv0\pmod5$, which, together with Ramanujan's
congruence $p(5n+4)\equiv0\pmod5$, implies that both the sum
of all odd parts and the sum of all even parts in the partitions
of $5n+4$ are divisible by $5$. They also conjectured that
$SOME(\lambda)\equiv0\pmod{5^\alpha}$, where $\alpha$ is a
positive integer and $\lambda$ is a nonnegative integer
satisfying $24\lambda\equiv1\pmod{5^\alpha}$, which was subsequently proved by Bardhan and Saikia \cite{BS}. Furthermore, Garvan and Sarma \cite{GarvanSarma2026} studied sums of certain parts in overpartitions and partitions without repeated odd parts. Subsequently, Baruah and Gogoi \cite{BaruahGogoi2026} investigated analogous sums of overlined parts in overpartitions. 

Motivated by the $SOME(n)$ function and its properties, in this paper, we define the corresponding arithmetic functions $COS(n)$ and $CES(n)$ for cubic partitions and study their divisibility properties. The function $COS(n)$ (resp. $CES(n)$) denotes the sum of all odd (resp. even) parts in all cubic partitions of a positive integer $n$. Since the sum of all parts of each cubic partition of $n$ is $n$, it follows that
\begin{equation}\label{e4}
	COS(n)+CES(n)=n a(n).
\end{equation}
For example, $COS(3)=8$ and $CES(3)=4$. For convenience, we set $COS(0)=CES(0)=0$. 

The paper is organized as follows. Section~2 records the necessary $q$-identities and extraction formulas used in the subsequent proofs. In Section~3, we derive the generating functions and divisor-sum formulas for $COS(n)$ and $CES(n)$ and establish their congruences modulo $3$, $4$, $5$, and $7$. Finally, Section~4 presents a conjecture for $COS(n)$ and $CES(n)$ modulo higher prime powers. 

\section{Preliminary results}

Ramanujan's general theta function is defined by \cite[Eq.~(1.2.1), p.~6]{Berndt2006}
$$
f(\alpha,\beta)=\sum_{m\in\Z}
\alpha^{m(m+1)/2}\beta^{m(m-1)/2},\qquad |\alpha\beta|<1.
$$
Three important special cases \cite[Eqs.~(1.2.2)--(1.2.4), pp.~6--7]{Berndt2006} are
\begin{align}\label{e5}
\phi(q)&:=f(q,q)
=\sum_{m\in\Z}q^{m^2}
=\dfrac{f_2^5}{f_1^2f_4^2},\\
\psi(q)&\label{e6}:=f(q,q^3)
=\sum_{m\ge0}q^{m(m+1)/2}
=\dfrac{f_2^2}{f_1},\\
\intertext{and }
f(-q)&\label{e7}:=f(-q,-q^2)
=\sum_{m\in\Z}(-1)^mq^{m(3m-1)/2}
=f_1.
\end{align}
The product representations in \eqref{e5}, \eqref{e6}, and \eqref{e7} follow from Jacobi's triple product identity \cite[Theorem~1.3.3]{Berndt2006}, given by
$$
f(\alpha,\beta)=(-\alpha;\alpha\beta)_\infty
(-\beta;\alpha\beta)_\infty(\alpha\beta;\alpha\beta)_\infty, 
$$
which can be equivalently stated as
\begin{equation}\label{e8}
\sum_{m\in\Z}z^mt^{m(m-1)/2}
=(-z,-t/z;t)_\infty(t;t)_\infty,
\qquad 0<|t|<1,\quad z\ne0.
\end{equation}
Replacing $q$ by $-q$ in \eqref{e5} and simplifying, we obtain
\begin{equation}\label{e9}
\phi(-q)=\sum_{m\in\Z}(-1)^mq^{m^2}=\dfrac{f_1^2}{f_2}.
\end{equation}
Replacing $q$ by $q^2$ in \eqref{e6}, we obtain
\begin{equation}\label{e10}
\psi(q^2)=\dfrac{f_4^2}{f_2}.
\end{equation}
Also, the binomial theorem gives, for every prime $p$,
\begin{equation}\label{e14}
	f_r^p=\prod_{m\ge1}(1-q^{rm})^p\equiv\prod_{m\ge1}(1-q^{prm})=f_{pr}\pmod p.
\end{equation}
The following identity due to Jacobi \cite[Theorem~1.3.9]{Berndt2006} will also be used in the proofs:
\begin{equation}\label{e11}
f_1^3=\sum_{m\ge0}(-1)^m(2m+1)q^{m(m+1)/2}.
\end{equation}

Next, we record some results related to Lambert series. Define 
\begin{equation}\label{e12}
L(q)=\sum_{m\ge1}\dfrac{mq^m}{1-q^m}
=\sum_{n\ge1}\sigma(n)q^n, 
\end{equation}
where $\sigma(n)$ denotes the sum of all positive divisors of the positive integer $n$. 
Throughout the paper, we will use the notation $D:=q\dfrac{d}{dq}$. By logarithmic differentiation of the $q$-product, $f_r=(q^r;q^r)_\infty$, we obtain
\begin{equation}\label{e13}
	\dfrac{Df_r}{f_r}=-r\sum_{m\ge1}\dfrac{mq^{rm}}{1-q^{rm}}=-rL(q^r).
\end{equation}

Throughout the rest of the paper, we will denote the coefficient of $q^n$ in $F(q)$ by $[q^n]F(q)$, where $F(q)$ is a formal power series in the variable $q$ with coefficients from the ring of integers $\Z$. We will use the standard notation $\Z[[q]]$ to denote the  ring of formal power series in the variable $q$ with coefficients from  $\Z$. 
All congruences between formal power series are understood coefficientwise. The operator $D$ preserves these congruences: if $F(q)-G(q)=mC(q)$ for a positive integer $m$ and some $C(q)\in\Z[[q]]$, then $DF(q)-DG(q)=mDC(q)\in m\Z[[q]]$.
Also, for any positive integer $m$ and $0\le r<m$, let the extraction operator $\mathcal E_{m,r}$ be defined by
   \begin{equation}\label{e15}
   	\mathcal E_{m,r}\left(\sum_{n\ge0}g(n)q^n\right)
   	=\sum_{n\ge0}g(mn+r)q^n.
   \end{equation}
Let $B(q)=\sum_{j\ge0}b(j)q^j$ and $G(q)=\sum_{n\ge0}g(n)q^n$.  Then for any prime $p$, 
\begin{align}\label{e16}
&\mathcal E_{p,r}\!\left(B(q^p)G(q)\right)=\sum_{n\ge0}\sum_{j=0}^nb(j)g(p(n-j)+r)q^n=B(q)\mathcal E_{p,r}(G(q))\\\intertext{and}
&\label{e17}\mathcal E_{p,r}\!\left((D-r)G(q)\right)=\sum_{n\ge0}pn\,g(pn+r)q^n=pD\bigl(\mathcal E_{p,r}(G(q))\bigr).
\end{align}
Further, if $G(q)\in\Z[[q]]$, then  from \eqref{e17}, we note that
\begin{equation}\label{e18}
\mathcal E_{p,r}\!\left((D-r)G(q)\right)\equiv0\pmod p
\end{equation}
and for $B(q)\in\Z[[q]]$, from the  binomial theorem and Fermat's little theorem, we obtain
\begin{equation}\label{e19}
B(q)^p\equiv\sum_{j\ge0}b(j)^pq^{pj}\equiv\sum_{j\ge0}b(j)q^{pj}=B(q^p)\pmod p.
\end{equation}
Moreover,  it can be verified that
\begin{equation}\label{e20}
\mathcal E_{p,s}\bigl(\mathcal E_{p,r}(G(q))\bigr)
=\sum_{n\ge0}g(p^2n+ps+r)q^n,
\qquad 0\le r,s<p.
\end{equation}

Next, we record some lemmas involving the Lambert series $L(q)$ and the three important Eisenstein series $E_2(q)$, $E_4(q)$, and $E_6(q)$ \cite{RSE}, which are defined by 
\begin{align}\label{e21}
E_2(q)&=1-24\sum_{n\ge1}\sigma_1(n)q^n=1-24L(q),\\
E_4(q)&\label{e22}=1+240\sum_{n\ge1}\sigma_3(n)q^n,\\\intertext{and}
E_6(q)&\label{e23}=1-504\sum_{n\ge1}\sigma_5(n)q^n, 
\end{align}
where, for any positive integer $k$, $\sigma_k(n)=\sum_{d\mid n}d^k$. Note that $\sigma_1(n)=\sigma(n)$.
The above three Eisenstein series satisfy the following differential equations due to Ramanujan \cite[Eq.~(30), \S~7]{RSE}:
\begin{align}\label{e24}
DE_2(q)&=\dfrac{E_2(q)^2-E_4(q)}{12},\\
DE_4(q)&\label{e25}=\dfrac{E_2(q)E_4(q)-E_6(q)}{3},\\\intertext{and}
DE_6(q)&\label{e26}=\dfrac{E_2(q)E_6(q)-E_4(q)^2}{2}.
\end{align}
From \eqref{e21}, it follows that
\begin{equation}\label{e27}
M:=M(q):=2E_2(q^2)-E_2(q)=1+24\bigl(L(q)-2L(q^2)\bigr).
\end{equation}

\begin{lemma}
We have
\begin{align}\label{e28}
&M=2\phi(q)^4-\phi(-q)^4=\phi(-q)^4+32q\psi(q^2)^4,\\
&\label{e29}5M^2=E_4(q)+4E_4(q^2),\\
&\label{e30}7M^3=8E_6(q^2)-E_6(q),\\\intertext{and}
&\label{e31}E_6(q^2)=M^3-72qM\dfrac{f_2^{16}}{f_1^8}.
\end{align}
\end{lemma}

\begin{proof}
Jacobi's four-square theorem \cite[p.~188]{Jacobi1829} states that the number of ordered representations of a positive integer $n$ as a sum of four integer squares, with signs and zero coordinates included, is $8\sum_{d\mid n,\,4\nmid d}d$. Thus, by \eqref{e12}, we obtain
\begin{align}
\phi(q)^4
&=1+8\sum_{n\ge1}\left(\sum_{\substack{d\mid n\\4\nmid d}}d\right)q^n\notag\\
&=1+8\sum_{n\ge1}\sigma(n)q^n
-8\sum_{n\ge1}\sum_{\substack{d\mid n\\4\mid d}}dq^n\notag\\
&=1+8L(q)-32\sum_{m\ge1}\left(\sum_{e\mid m}e\right)q^{4m}\notag\\
&\label{e32}=1+8L(q)-32L(q^4).
\end{align}
Set $n=2^au$, where $a\ge0$ and $u$ is an odd integer. Then 
\begin{equation}\label{e33}
\sigma(2^au)=\sum_{j=0}^{a}\sum_{d\mid u}2^jd=(2^{a+1}-1)\sigma(u).
\end{equation}
From \eqref{e33}, it follows that, for $a=1$,  $$\sigma(2u)=3\sigma(u)$$ and for $a\ge1$, 
$$
3\sigma(2^au)-2\sigma(2^{a-1}u)=\bigl(3(2^{a+1}-1)-2(2^a-1)\bigr)\sigma(u)
=\sigma(2^{a+1}u).
$$
Therefore,  for every integer $n\ge1$, we obtain
\begin{equation}\label{e34}
\sigma(2n)=3\sigma(n)-2\sigma(n/2), 
\end{equation}
where we used the convention that $\sigma(x)=0$ when $x$ is not a positive integer.
Thus,  from \eqref{e12} and \eqref{e34}, we obtain
\begin{align}
L(-q)
&=\sum_{n\ge1}(-1)^n\sigma(n)q^n\notag\\
&=-\sum_{n\ge1}\sigma(n)q^n
+2\sum_{\substack{n\ge1\\n\ {\rm even}}}\sigma(n)q^n\notag\\
&=-L(q)+2\sum_{m\ge1}\sigma(2m)q^{2m}\notag\\
&=-L(q)+6\sum_{m\ge1}\sigma(m)q^{2m}
-4\sum_{m\ge1}\sigma(m/2)q^{2m}\notag\\
&=-L(q)+6L(q^2)-4\sum_{r\ge1}\sigma(r)q^{4r}\notag\\
&\label{e35}=-L(q)+6L(q^2)-4L(q^4).
\end{align}
Replacing $q$ by $-q$ in \eqref{e32} and then employing \eqref{e35}, we obtain
\begin{equation}\label{e36}
\phi(-q)^4=1+8L(-q)-32L(q^4)=1-8L(q)+48L(q^2)-64L(q^4).
\end{equation}
From \eqref{e27}, \eqref{e32}, and \eqref{e36}, we obtain
\begin{equation}\label{p1}
2\phi(q)^4-\phi(-q)^4
=1+24L(q)-48L(q^2)
=M.
\end{equation}
From \cite[Eq.~(3.6.8), p.~72]{Berndt2006}, we note that
\begin{equation}\label{e37}
\phi(q)^4-\phi(-q)^4=16q\psi(q^2)^4.
\end{equation}
Employing \eqref{e37} in \eqref{p1}, we obtain
$$
M=\phi(-q)^4+2\bigl(\phi(q)^4-\phi(-q)^4\bigr)=\phi(-q)^4+32q\psi(q^2)^4,$$
which is \eqref{e28}.

To prove the remaining identities, set
\begin{equation}\label{e38}
z=\phi(q)^2 \qquad \text{and}\qquad 
x=\dfrac{\phi(q)^4-\phi(-q)^4}{\phi(q)^4}.
\end{equation}
Using \eqref{e28} and \eqref{e38}, we obtain
\begin{equation}\label{e39}
M=z^2(1+x)=\phi(q)^4+\phi(q)^4-\phi(-q)^4.
\end{equation}
From \cite[pp.~126--127, Theorems~5.4.11 and 5.4.12]{Berndt2006}, we note that
\begin{align}\label{e40}
&E_4(q)=z^4(1+14x+x^2),\\
&\label{e41}E_4(q^2)=z^4(1-x+x^2),\\
&\label{e42}E_6(q)=z^6(1+x)(1-34x+x^2),\\\intertext{and}
&\label{e43}E_6(q^2)=z^6(1+x)(1-x/2)(1-2x).
\end{align}
Combining \eqref{e39},  \eqref{e40}, and \eqref{e41}, we obtain
\begin{align*}
E_4(q)+4E_4(q^2)
=z^4\bigl(1+14x+x^2+4-4x+4x^2\bigr)
=5\bigl(z^2(1+x)\bigr)^2
=5M^2,
\end{align*}
which is \eqref{e29}.\\
Again, using \eqref{e39}, \eqref{e42}, and \eqref{e43}, we obtain
\begin{align*}
8E_6(q^2)-E_6(q)
=z^6(1+x)(7+14x+7x^2)=7\bigl(z^2(1+x)\bigr)^3=7M^3,
\end{align*}
which is \eqref{e30}.

Finally, using \eqref{e5}, \eqref{e10}, \eqref{e37}, \eqref{e38}, \eqref{e39}, and \eqref{e43}, we obtain
\begin{align*}
M^3-E_6(q^2)
&=z^6(1+x)\bigl((1+x)^2-(1-x/2)(1-2x)\bigr)\\
&=\dfrac92xz^6(1+x)\\
&=\dfrac92M\phi(q)^4\bigl(\phi(q)^4-\phi(-q)^4\bigr)\\
&=72qM\phi(q)^4\psi(q^2)^4\\
&=72qM\dfrac{f_2^{20}}{f_1^8f_4^8}\cdot\dfrac{f_4^8}{f_2^4}\\
&=72qM\dfrac{f_2^{16}}{f_1^8}. 
\end{align*}
This completes the proof of \eqref{e31}.
\end{proof}

\begin{lemma}
We have
\begin{equation}\label{e50}
L(q)\equiv L(q^3)+qf_3^8\pmod3
\end{equation}
and the following exact dissection holds:
\begin{equation}\label{e51}
\psi(q)=\dfrac{f_6f_9^2}{f_3f_{18}}+q\psi(q^9).
\end{equation}

\end{lemma}

\begin{proof}
In the series on the right-hand side of Jacobi's identity \eqref{e11}, the exponent $r(r+1)/2$ is congruent to $0$ modulo $3$ if $r\equiv0\text{ or }2\pmod3$ and is congruent to $1$ modulo $3$ if $r\equiv1\pmod3$. Thus, putting $r=3u+1$ with $u\ge0$, we obtain
\begin{equation}\label{e52}
\sum_{u\ge0}(-1)^{3u+1}(6u+3)q^{(3u+1)(3u+2)/2}=-3q\sum_{u\ge0}(-1)^u(2u+1)q^{9u(u+1)/2}=-3qf_9^3
\end{equation}
Therefore, from \eqref{e11} and \eqref{e52}, we obtain
\begin{equation}\label{e53}
f_1^3=B(q^3)-3qf_9^3, 
\end{equation}for some $B(q)\in\Z[[q]]$.
Using \eqref{e14} in \eqref{e53}, we note that $$B(q^3)\equiv f_3\pmod3,$$ which immediately implies
$$B(q)\equiv f_1\pmod3.
$$
Applying the differential operator $D$ to $B(q)$, replacing $q$ by $q^3$, and employing \eqref{e13}, we obtain
\begin{equation}\label{e54}
\left.DB(q)\right|_{q\mapsto q^3}\equiv\left.Df_1\right|_{q\mapsto q^3}=-L(q^3)f_3\pmod3.
\end{equation}
Differentiating \eqref{e53} and simplifying  by using \eqref{e13}, we obtain
\begin{equation}\label{e55}
-3L(q)f_1^3=3\left.DB(q)\right|_{q\mapsto q^3}-3qf_9^3-3qD(f_9^3).
\end{equation}
Dividing the exact identity \eqref{e55} by $-3$, we obtain
$$
L(q)f_1^3=-\left.DB(q)\right|_{q\mapsto q^3}+qf_9^3+qD(f_9^3).
$$
Since $D(f_9^3)=-27L(q^9)f_9^3$, reduction modulo $3$, followed by \eqref{e14} and \eqref{e54}, gives
$$
L(q)f_3\equiv L(q^3)f_3+qf_9^3\pmod3.
$$
Since $f_3$ has constant term $1$, we obtain
\begin{equation*}
L(q)\equiv L(q^3)+q\dfrac{f_9^3}{f_3}\equiv L(q^3)+q\dfrac{f_3^9}{f_3}=L(q^3)+qf_3^8\pmod3,
\end{equation*}
which is  \eqref{e50}.

Next,  express \eqref{e6} as
\begin{equation}\label{e56}
\psi(q)
=\sum_{u\ge0}q^{3u(3u+1)/2}
+q\sum_{u\ge0}q^{9u(u+1)/2}
+\sum_{u\ge0}q^{3(u+1)(3u+2)/2}=\sum_{t\in\Z}(q^3)^{t(3t+1)/2}+q\psi(q^9).
\end{equation}
Setting $t=q^3$ and $z=q^2$ in \eqref{e8}, we obtain
\begin{equation}\label{e57}
\sum_{n\in\Z}q^{n(3n+1)/2}=(-q^2,-q;q^3)_\infty f_3=\dfrac{f_2f_3^2}{f_1f_6}.
\end{equation}
Replacing $q$ by $q^3$ in \eqref{e57}, we obtain
\begin{equation}\label{nn1}
\sum_{t\in\Z}(q^3)^{t(3t+1)/2}=\dfrac{f_6f_9^2}{f_3f_{18}}.
\end{equation}
Employing \eqref{nn1} in \eqref{e56}, we complete the proof of  \eqref{e51}.
\end{proof}

\begin{lemma}
We have
\begin{align}\label{e81}
&\mathcal E_{5,2}(f_1^2)=-f_5^2,\\
&\label{e82}\mathcal E_{5,2}\!\left(\dfrac{f_1f_4^2}{f_2}\right)=q\dfrac{f_5f_{20}^2}{f_{10}},\\
&\label{e83}\mathcal E_{7,1}(f_2^2)=qf_{14}^2,\\\intertext{and}
&\label{e84}\mathcal E_{7,1}\!\left(\phi(q)f_4\right)=2\phi(-q)f_4+q\phi(q^7)f_{28}.
\end{align}
\end{lemma}

\begin{proof}
Set
$$
T(q)=\dfrac{(q^2,q^3;q^5)_\infty}{(q,q^4;q^5)_\infty}.
$$
From \cite[Chapter~16]{tifr} (see also \cite[Theorem~7.4.1, p.~161]{Berndt2006}), we note that
\begin{equation}\label{e85}
f_1=f_{25}\left(T(q^5)-q-\dfrac{q^2}{T(q^5)}\right).
\end{equation}
Squaring \eqref{e85}, we obtain
\begin{equation}\label{e86}
f_1^2
=f_{25}^2\left(T(q^5)^2-2qT(q^5)-q^2
+\dfrac{2q^3}{T(q^5)}+\dfrac{q^4}{T(q^5)^2}\right).
\end{equation}
Extracting the terms involving $q^{5n+2}$ from \eqref{e86}, dividing by $q^2$, and replacing $q^5$ by $q$, we arrive at \eqref{e81}.

Replacing $q$ by $q^4$ in \eqref{e7}, we obtain
\begin{equation}\label{e87}
f_4=\sum_{r\in\Z}(-1)^rq^{6r^2-2r}.\end{equation}
Setting $t=q^4$ and $z=-q$ in \eqref{e8}, we obtain
\begin{equation}
\label{e88}\sum_{s\in\Z}(-1)^sq^{2s^2-s}=(q,q^3;q^4)_\infty f_4=\dfrac{f_1f_4}{f_2}.
\end{equation}
Multiplying \eqref{e87} and \eqref{e88}, we obtain
\begin{equation}\label{e89}
\dfrac{f_1f_4^2}{f_2}
=\sum_{r,s\in\Z}(-1)^{r+s}q^{6r^2-2r+2s^2-s}.
\end{equation}
Set
$$
N_1=6r^2-2r+2s^2-s.
$$Then
\begin{equation}\label{e90}
24N_1+7=(12r-2)^2+3(4s-1)^2.
\end{equation}
Since $-3\equiv2\pmod5$ is a quadratic nonresidue modulo $5$, \eqref{e90} implies $
N_1\equiv2\pmod5$ if and only if 
$r=1-5u$ and 
$s=5v-1$,
 where $u$ and $v$ are integers.
 Consequently, from \eqref{e89}, we obtain
\begin{align*}
\mathcal E_{5,2}\!\left(\dfrac{f_1f_4^2}{f_2}\right)
&=\sum_{u,v\in\Z}(-1)^{u+v}
q^{(7+25(6u^2-2u+2v^2-v)-2)/5}\\
&=q\left(\sum_{u\in\Z}(-1)^uq^{5(6u^2-2u)}\right)
\left(\sum_{v\in\Z}(-1)^vq^{5(2v^2-v)}\right)\\
&=qf_{20}\dfrac{f_5f_{20}}{f_{10}}\\
&=q\dfrac{f_5f_{20}^2}{f_{10}},
\end{align*}
which completes the proof of \eqref{e82}.

Replacing $q$ by $q^2$ in \eqref{e7} and squaring, we obtain
\begin{equation}\label{e94}
f_2^2=\sum_{r,s\in\Z}(-1)^{r+s}q^{3r^2-r+3s^2-s}.
\end{equation}
Set
$$N_2=3r^2-r+3s^2-s.$$
Then
\begin{equation}\label{e95}
12N_2+2=(6r-1)^2+(6s-1)^2.
\end{equation}
Since $-1$ is a quadratic nonresidue modulo $7$, \eqref{e95} gives
$
N_2\equiv1\pmod7$ if and only if 
$r=7u-1$ and 
$s=7v-1$.
Therefore, 
$$
\mathcal E_{7,1}(f_2^2)=\sum_{u,v\in\Z}(-1)^{u+v}
q^{(8+49(3u^2-u+3v^2-v)-1)/7}=q\left(\sum_{u\in\Z}(-1)^uq^{7(3u^2-u)}\right)^2=qf_{14}^2,
$$
which completes the proof of \eqref{e83}.

From \eqref{e5} and \eqref{e87}, we obtain
\begin{equation}\label{e99}
\phi(q)f_4=\sum_{r,s\in\Z}(-1)^sq^{r^2+6s^2-2s}.
\end{equation}
The condition $r^2+6s^2-2s\equiv1\pmod7$ is equivalent to $r\equiv\pm(s+1)\pmod7$. The class $r\equiv-(s+1)\pmod7$ is parametrized bijectively by
$$
r=6t-v-1,\qquad s=t+v,\qquad t,v\in\Z,
$$
since $t=(r+s+1)/7$ and $v=s-t$. Replacing $r$ by $-r$ parametrizes the other class. The two classes intersect when $r\equiv0\pmod7$ and $s\equiv-1\pmod7$, so this intersection must be counted only once. For either parametrization, we have
$$
r^2+6s^2-2s=1+7(v^2+6t^2-2t)\quad \text{and}\quad
(-1)^s=(-1)^{t+v}.
$$
So by \eqref{e9} and \eqref{e87}, we obtain
\begin{align*}
\sum_{t,v\in\Z}(-1)^{t+v}q^{v^2+6t^2-2t}
&=\left(\sum_{v\in\Z}(-1)^vq^{v^2}\right)
\left(\sum_{t\in\Z}(-1)^tq^{6t^2-2t}\right)=\phi(-q)f_4.
\end{align*}
Also,  when $r\equiv0\pmod7$ and $s\equiv-1\pmod7$,
writing $r=7d$ and $s=7u-1$ for $d,u\in\Z$, we obtain
$$
r^2+6s^2-2s=8+49(d^2+6u^2-2u)\quad \text{and}\quad
(-1)^s=-(-1)^u.
$$
So, when $r\equiv0\pmod7$ and $s\equiv-1\pmod7$, 
\begin{align*}
-\sum_{d,u\in\Z}(-1)^u q^{1+7(d^2+6u^2-2u)}
&=-q\phi(q^7)f_{28}.
\end{align*}
Therefore, we obtain
\begin{align*}
\mathcal E_{7,1}\!\left(\phi(q)f_4\right)
&=2\phi(-q)f_4-\bigl(-q\phi(q^7)f_{28}\bigr)=2\phi(-q)f_4+q\phi(q^7)f_{28}. 
\end{align*}
Hence, the proof of  \eqref{e84} is complete.
\end{proof}

\begin{lemma}
We have
\begin{align}\label{e112}
&\mathcal E_{5,2}(M A(q))\equiv Mf_1^3f_2^3\pmod5\\\intertext{and}
&\label{e113}\mathcal E_{7,1}(MA(q))\equiv4f_1f_2M^3\pmod7.
\end{align}
\end{lemma}

\begin{proof}
Set 
$$K=\dfrac{f_2^{19}}{f_1^9f_4^8}\quad \text{and}\quad J=\dfrac{f_1^7}{f_2^5}.$$
From \eqref{e2}, \eqref{e5}, \eqref{e9}, and \eqref{e28}, we obtain
\begin{equation}\label{e114}
MA(q)=\dfrac{2\phi(q)^4-\phi(-q)^4}{f_1f_2}=2K-J.
\end{equation}
Using \eqref{e14}, we obtain
\begin{align}\label{e115}
&K=\dfrac{f_2^{20}}{f_1^{10}f_4^{10}}\dfrac{f_1f_4^2}{f_2}\equiv\dfrac{f_{10}^4}{f_5^2f_{20}^2}\dfrac{f_1f_4^2}{f_2}\pmod5,\\\intertext{and}
&\label{e116}J=\dfrac{f_1^5}{f_2^5}f_1^2\equiv\dfrac{f_5}{f_{10}}f_1^2\pmod5.
\end{align}
Employing \eqref{e115} and \eqref{e116} in \eqref{e114}, and then simplifying by using \eqref{e16}, \eqref{e81}, and \eqref{e82}, we obtain
\begin{equation}\label{e117}
\mathcal E_{5,2}(MA(q))
\equiv2\dfrac{f_2^4}{f_1^2f_4^2}
\left(q\dfrac{f_5f_{20}^2}{f_{10}}\right)
-\dfrac{f_1}{f_2}(-f_5^2)=2q\dfrac{f_2^4f_5f_{20}^2}{f_1^2f_4^2f_{10}}
+\dfrac{f_1f_5^2}{f_2}\pmod5.
\end{equation}
Employing \eqref{e9}, \eqref{e10}, \eqref{e14}, and \eqref{e28} in \eqref{e117}, we obtain
\begin{align*}
\mathcal E_{5,2}(MA(q))
&\equiv2q\dfrac{f_2^4}{f_1^2f_4^2}
\dfrac{f_1^5f_4^{10}}{f_2^5}
+\dfrac{f_1}{f_2}f_1^{10}\pmod5\\
&=2q\dfrac{f_1^3f_4^8}{f_2}+\dfrac{f_1^{11}}{f_2}\\
&=f_1^3f_2^3\left(2q\dfrac{f_4^8}{f_2^4}
+\dfrac{f_1^8}{f_2^4}\right)\\
&=f_1^3f_2^3\bigl(2q\psi(q^2)^4+\phi(-q)^4\bigr)\\
&\equiv Mf_1^3f_2^3\pmod5.
\end{align*}
So the proof of \eqref{e112} is complete.

Similarly, from \eqref{e14} and \eqref{e5}, we obtain
\begin{align}\label{e118}
K&=\dfrac{f_2^{14}}{f_1^7f_4^7}\dfrac{f_2^5}{f_1^2f_4}\equiv\dfrac{f_{14}^2}{f_7f_{28}}\phi(q)f_4\pmod7,\\\intertext{and}
J&\label{e119}=\dfrac{f_1^7}{f_2^7}f_2^2\equiv\dfrac{f_7}{f_{14}}f_2^2\pmod7.
\end{align}
Employing \eqref{e16} and \eqref{e84} in \eqref{e118}, and then using \eqref{e5}, \eqref{e9}, \eqref{e10}, \eqref{e14}, and \eqref{e19}, we obtain
\begin{align}
\mathcal E_{7,1}(K)&\equiv\dfrac{f_2^2}{f_1f_4}
\bigl(2\phi(-q)f_4+q\phi(q^7)f_{28}\bigr)\pmod7\notag\\
&\equiv2f_1f_2+q\dfrac{f_2^2f_4^6}{f_1}\phi(q)^7\pmod7\notag\\
&=2f_1f_2+q\dfrac{f_2^{37}}{f_1^{15}f_4^8}\notag\\
&=f_1f_2\left(2+q\dfrac{f_4^8}{f_2^4}
\dfrac{f_2^{40}}{f_1^{16}f_4^{16}}\right)\notag\\
&\label{e120}=f_1f_2\bigl(2+q\psi(q^2)^4\phi(q)^8\bigr).
\end{align}
Similarly, employing \eqref{e16} and \eqref{e83} in \eqref{e119}, and then using \eqref{e14}, \eqref{e5}, \eqref{e9}, and \eqref{e10}, we obtain
\begin{align}
\mathcal E_{7,1}(J)
&\equiv\dfrac{f_1}{f_2}\,qf_{14}^2\pmod7\notag\\
&\equiv qf_1f_2^{13}\pmod7\notag\\
&=qf_1f_2\,\dfrac{f_4^8}{f_2^4}
\cdot\dfrac{f_2^{20}}{f_1^8f_4^8}\cdot\dfrac{f_1^8}{f_2^4}\notag\\
&\label{e121}=qf_1f_2\psi(q^2)^4\phi(q)^4\phi(-q)^4.
\end{align}
Employing \eqref{e120} and \eqref{e121} in \eqref{e114}, and then using \eqref{e28}, \eqref{e5}, and \eqref{e10}, we obtain
\begin{align}
\mathcal E_{7,1}(M\cdot A(q))
&\equiv f_1f_2\bigl(4+2q\psi(q^2)^4\phi(q)^8
-q\psi(q^2)^4\phi(q)^4\phi(-q)^4\bigr)\pmod7\notag\\
&=f_1f_2\bigl(4+q\psi(q^2)^4\phi(q)^4
(2\phi(q)^4-\phi(-q)^4)\bigr)\notag\\
&\label{e122}=f_1f_2\left(4+qM\dfrac{f_2^{16}}{f_1^8}\right).
\end{align}
From \eqref{e23}, it is easily seen that  $E_6(q^2)\equiv1\pmod7$. So  from \eqref{e31}, we obtain
\begin{equation}\label{e123}
qM\dfrac{f_2^{16}}{f_1^8}\equiv4(M^3-1)\pmod7.
\end{equation}
Employing \eqref{e123} in \eqref{e122}, we obtain
$$\mathcal E_{7,1}(MA(q))\equiv f_1f_2\bigl(4+4(M^3-1)\bigr)=4f_1f_2M^3\pmod7,$$
which completes the proof of  \eqref{e113}.
\end{proof}

\begin{lemma}
We have
\begin{align}\label{e124}
&\mathcal E_{5,2}(A(q))\equiv3f_1^3f_2^3\pmod5\\\intertext{and}
&\label{e125}\mathcal E_{7,1}(A(q))\equiv f_1f_2M^2\pmod7.
\end{align}
\end{lemma}

\begin{proof}
First, we show that
\begin{align}\label{e126}
(D-2)(M^2A(q))&\equiv(M^5+2M)A(q)\pmod5\\\intertext{and}
(D-1)(M^3A(q))&\label{e127}\equiv(M^7+5M)A(q)\pmod7.
\end{align}
Set
$$P=E_2(q),\quad Q=E_2(q^2),\quad S=E_4(q),\quad R=E_4(q^2).$$
From \eqref{e2} and \eqref{e13}, we obtain
\begin{equation}\label{n3}
\dfrac{DA(q)}{A(q)}=-\dfrac{Df_1}{f_1}-\dfrac{Df_2}{f_2}=L(q)+2L(q^2).
\end{equation}
Employing \eqref{e21} in \eqref{n3}, we obtain
\begin{equation}\label{e128}
\dfrac{DA(q)}{A(q)}=\dfrac{1-P}{24}+2\dfrac{1-Q}{24}=\dfrac{3-P-2Q}{24}.
\end{equation}
Differentiating \eqref{e27} and using \eqref{e24}, we obtain
\begin{align}
DM&=2D\bigl(E_2(q^2)\bigr)-DE_2(q)\notag\\
&=\dfrac{4(Q^2-R)-(P^2-S)}{12}\notag\\
&=\dfrac{(2Q-P)(2Q+P)+S-4R}{12}\notag\\
&\label{e129}=\dfrac{M(P+2Q)+S-4R}{12}.
\end{align}
From \eqref{e22}, it is easily seen that 
\begin{equation}\label{sv}S\equiv 1\pmod5 \quad \text{and}\quad  R\equiv1\pmod5.\end{equation} So from  \eqref{e25}, we obtain
\begin{align*}
0&\equiv3DS=PS-E_6(q)\equiv P-E_6(q)\pmod5\\\intertext{and}
0&\equiv3DR=2QR-2E_6(q^2)\equiv2\bigl(Q-E_6(q^2)\bigr)\pmod5,
\end{align*}which imply, respectively,
\begin{align}\label{e130}
&E_6(q)\equiv P\pmod5\\\intertext{and}
&\label{e131}E_6(q^2)\equiv Q\pmod5.
\end{align}
Employing \eqref{e130} and \eqref{e131} in \eqref{e30}, and then using \eqref{e27}, we obtain
\begin{align}
2M^3&\equiv3Q-P=Q+M\pmod5,\notag\\
Q&\label{e132}\equiv2M^3-M\pmod5,\\\intertext{and}
P&\label{e133}=2Q-M\equiv4M^3-3M\pmod5.
\end{align}
Using \eqref{e132} and \eqref{e133}, we obtain
\begin{equation}\label{e134}
P+2Q\equiv(4M^3-3M)+2(2M^3-M)=8M^3-5M\equiv3M^3\pmod5.
\end{equation}
Employing \eqref{e134} in \eqref{e128} and \eqref{e129}, and noting $S\equiv R\equiv1\pmod5$, we obtain
\begin{align}
&\label{e136}\dfrac{DA(q)}{A(q)}\equiv4(3-3M^3)\equiv2+3M^3\pmod5\\\intertext{and}
&\label{e135}DM\equiv3\bigl(3M^4+1-4\bigr)\equiv1-M^4\pmod5,
\end{align}respectively.
Combining \eqref{e136} and \eqref{e135}, we obtain
\begin{equation*}
\dfrac{(D-2)(M^2A(q))}{A(q)}
=2M\,DM+M^2\dfrac{DA(q)}{A(q)}-2M^2
\equiv M^5+2M\pmod 5,
\end{equation*}
which is equivalent to \eqref{e126}.

Again, from  \eqref{e23}, we note that
\begin{equation*}
E_6(q)\equiv1\pmod7\quad \text{and}\quad
E_6(q^2)\equiv1\pmod7.
\end{equation*}
Therefore, from \eqref{e26}, we obtain
\begin{align*}
0&\equiv2DE_6(q)=PE_6(q)-S^2\equiv P-S^2\pmod7\\\intertext{and}
0&\equiv D\bigl(E_6(q^2)\bigr)=QE_6(q^2)-R^2\equiv Q-R^2\pmod7,
\end{align*}
which imply
\begin{align}\label{e137}
&P\equiv S^2\pmod7\\\intertext{and}
&\label{e138}Q\equiv R^2\pmod7,
\end{align}respectively.
Employing \eqref{e137} and \eqref{e138} in \eqref{e27}, and then using \eqref{e29}, we obtain
\begin{equation}\label{e139}
M\equiv2R^2-S^2=2R^2-(5M^2-4R)^2
\equiv3M^4+5M^2R\pmod7.
\end{equation}
Since $M$ has constant term $1$, it is a unit in $\Z[[q]]$. Hence its inverse exists in $\Z[[q]]$ and remains an inverse after reduction modulo $5$ or $7$. Multiplying \eqref{e139} by $3M^{-2}$, and then using \eqref{e29}, we obtain
\begin{align}\label{e140}
R&\equiv3M^{-2}(M-3M^4)\equiv5M^2+3M^{-1}\pmod7\\\intertext{and}
S&\label{e141}=5M^2-4R\equiv5M^2-4(5M^2+3M^{-1})\equiv6M^2+2M^{-1}\pmod7.
\end{align}
Employing \eqref{e140} and \eqref{e141} in \eqref{e137} and \eqref{e138}, we obtain
\begin{align}\label{e142}
&P\equiv(6M^2+2M^{-1})^2=36M^4+24M+4M^{-2}\equiv M^4+3M+4M^{-2}\pmod7\\\intertext{and}
&\label{e143}Q\equiv(5M^2+3M^{-1})^2=25M^4+30M+9M^{-2}\equiv4M^4+2M+2M^{-2}\pmod7,
\end{align}respectively. From \eqref{e140} and \eqref{e141}, we obtain
\begin{equation}\label{e145}S-4R\equiv-14M^2-10M^{-1}\equiv4M^{-1}\pmod7.
	\end{equation}
From  \eqref{e142} and \eqref{e143}, we obtain
\begin{equation}\label{e144}
P+2Q\equiv9M^4+7M+8M^{-2}\equiv2M^4+M^{-2}\pmod7.
\end{equation}
Employing \eqref{e145} and \eqref{e144} in \eqref{e128} and \eqref{e129}, we obtain
\begin{align}
&\label{e147}\dfrac{DA(q)}{A(q)}\equiv5(3-2M^4-M^{-2})\equiv1+4M^4+2M^{-2}\pmod7\\\intertext{and}
\label{e146}
&DM\equiv3\bigl(M(2M^4+M^{-2})+4M^{-1}\bigr)=6M^5+15M^{-1}\equiv M^{-1}-M^5\pmod7.
\end{align}respectively.
Combining \eqref{e147} and \eqref{e146}, we obtain
\begin{align*}
\dfrac{(D-1)(M^3A(q))}{A(q)}
&=3M^2\,DM+M^3\left(\dfrac{DA(q)}{A(q)}\right)-M^3\\
&\equiv3M^2(M^{-1}-M^5)+M^3(1+4M^4+2M^{-2})-M^3\\
&=M^7+5M\pmod 7,
\end{align*}
which is equivalent to \eqref{e127}.

Now, 
employing \eqref{e18}, \eqref{e19}, \eqref{e16}, and \eqref{e112} in \eqref{e126}, we obtain
\begin{align}
0&\equiv\mathcal E_{5,2}\bigl((D-2)(M^2A(q))\bigr)\notag\\
&\equiv\mathcal E_{5,2}(M^5A(q))+2\mathcal E_{5,2}(MA(q))\notag\\
&\equiv\mathcal E_{5,2}\bigl(M(q^5)A(q)\bigr)+2Mf_1^3f_2^3\notag\\
&\label{e148}=M\mathcal E_{5,2}(A(q))+2Mf_1^3f_2^3\pmod 5.
\end{align}
Multiplying \eqref{e148} by $M^{-1}$, we arrive at
 \eqref{e124}.

Similarly, employing \eqref{e18}, \eqref{e19}, \eqref{e16}, and \eqref{e113} in \eqref{e127}, we obtain
\begin{align}
0&\equiv\mathcal E_{7,1}\bigl((D-1)(M^3A(q))\bigr)\notag\\
&\equiv\mathcal E_{7,1}(M^7A(q))+5\mathcal E_{7,1}(MA(q))\notag\\
&\equiv\mathcal E_{7,1}\bigl(M(q^7)A(q)\bigr)+20f_1f_2M^3\notag\\
&\label{e149}=M\mathcal E_{7,1}(A(q))+20f_1f_2M^3\pmod 7.
\end{align}
Multiplying \eqref{e149} by $M^{-1}$, we complete the proof of \eqref{e125}.
\end{proof}

\begin{lemma}
For $j=2,3$, let
\begin{equation}\label{e173}
H_j(q)=f_1f_2M^j.
\end{equation}
Then
\begin{equation}\label{e174}
\mathcal E_{7,6}(H_j(q))=7^{2j}f_7f_{14}M(q^7)^j.
\end{equation}
\end{lemma}

\begin{proof}
For $r\in\Z$, set
\begin{equation}\label{e175}
P_r=\dfrac{r(3r-1)}2.
\end{equation}
Ramanujan \cite[p.~369]{lost} recorded the following sequence of $q$-series, where $k$ is a nonnegative integer:

\begin{equation}\label{v2k}V_{2k}(q)=\dfrac{\sum_{r\in \Z} (-1)^r(6r-1)^{2k}q^{r(3r-1)/2}}{\sum_{r \in\Z}(-1)^rq^{r(3r-1)/2}}.\end{equation}
Using \eqref{e7} in \eqref{v2k}, we obtain
\begin{equation}\label{e176}
\sum_{r\in\Z}
(-1)^r(6r-1)^{2k}q^{P_r}
=
f_1V_{2k}(q).
\end{equation}
Since $(6r-1)^2=24P_r+1$, \eqref{e176} gives
\begin{align*}
f_1V_{2k+2}(q)&=(24D+1)\bigl(f_1V_{2k}(q)\bigr),\\
V_{2k+2}(q)&=24DV_{2k}(q)+E_2(q)V_{2k}(q),
\end{align*}
where we used \eqref{e13} and \eqref{e21}. Starting with $V_0(q)=1$ and employing \eqref{e24}--\eqref{e26}, we obtain the following identities recorded by Ramanujan \cite{lost}:
\begin{align}
V_0(q)&\label{e177a1} =1,\\
V_2(q)&\label{e177a2}=E_2(q),\\
V_4(q)&\label{e177a3} =3E_2(q)^2-2E_4(q),\\\intertext{and}
V_6(q)&\label{e177}=15E_2(q)^3-30E_2(q)E_4(q)+16E_6(q).
\end{align}
Now, for $r,s\in\Z$, set
\begin{equation}\label{e178}
x=6r-1,\quad y=6s-1,\quad\text{and}\quad  N=P_r+2P_s.
\end{equation}
Define
\begin{align}\label{e179}
W_2(x,y)&=x^4-12x^2y^2+4y^4,\\
W_3(x,y)&\label{e180}=x^6-30x^4y^2+60x^2y^4-8y^6,
\end{align}
and set
\begin{equation}\label{e181}
S_j(q)=\sum_{r,s\in\Z}(-1)^{r+s}W_j(x,y)q^N,
\qquad j=2,3.
\end{equation}
First we establish that  
\begin{align}\label{e182}
S_2(q)&=-7f_1f_2M^2\\\intertext{and}
S_3(q)&\label{e183}=23f_1f_2M^3.
\end{align}
Set $P=E_2(q)$ and $Q=E_2(q^2)$.

Employing \eqref{e176}--\eqref{e177}, and \eqref{e179} in \eqref{e181}, and then using \eqref{e27} and \eqref{e29}, we obtain
\begin{align*}
\dfrac{S_2(q)}{f_1f_2}
&=V_4(q)-12V_2(q)V_2(q^2)+4V_4(q^2)\\
&=3P^2-2E_4(q)-12PQ+12Q^2-8E_4(q^2)\\
&=3(P-2Q)^2-2\bigl(E_4(q)+4E_4(q^2)\bigr)\\
&=3M^2-2(5M^2)\\
&=-7M^2,
\end{align*}
which is equivalent to  \eqref{e182}.

Similarly, employing \eqref{e176}--\eqref{e177}, and \eqref{e180} in \eqref{e181}, and then using \eqref{e27}, \eqref{e29}, and \eqref{e30}, we obtain
\begin{align*}
\dfrac{S_3(q)}{f_1f_2}
&=V_6(q)-30V_4(q)Q+60PV_4(q^2)-8V_6(q^2)\\
&=15P^3-30PE_4(q)+16E_6(q)\\
&\quad-90P^2Q+60QE_4(q)+180PQ^2-120PE_4(q^2)\\
&\quad-120Q^3+240QE_4(q^2)-128E_6(q^2)\\
&=15(P-2Q)^3+30(2Q-P)E_4(q)\\
&\quad+120(2Q-P)E_4(q^2)+16\bigl(E_6(q)-8E_6(q^2)\bigr)\\
&=-15M^3+30M\bigl(E_4(q)+4E_4(q^2)\bigr)\\
&\quad+16\bigl(E_6(q)-8E_6(q^2)\bigr)\\
&=-15M^3+30M(5M^2)+16(-7M^3)\\
&=23M^3,
\end{align*}
which is equivalent to \eqref{e183}.

Now, from \eqref{e175} and \eqref{e178}, we obtain
\begin{equation}\label{e184}
x^2+2y^2=(24P_r+1)+2(24P_s+1)=24N+3.
\end{equation}
Since $-2\equiv5\pmod7$ is a quadratic nonresidue modulo $7$, from \eqref{e184}, we obtain
\begin{align}\label{e185}
N\equiv6\pmod7
&\Longleftrightarrow
\begin{cases}
x\equiv0\pmod7,\\
y\equiv0\pmod7,
\end{cases}\Longleftrightarrow
\begin{cases}
r=7u-1,\\
s=7v-1,
\end{cases}
\qquad u,v\in\Z.
\end{align}
Employing \eqref{e185} in \eqref{e175} and \eqref{e178}, we obtain
\begin{align}\label{e186}
P_{7u-1}&=49P_u+2,\\
\label{e187}N&=49(P_u+2P_v)+6,\\
\label{e188}(-1)^{r+s}&=(-1)^{7u+7v-2}=(-1)^{u+v},\\\intertext{and}
\label{e189}x&=7(6u-1), \quad y=7(6v-1).
\end{align}
Therefore,  employing \eqref{e185}, \eqref{e186}, \eqref{e187}, \eqref{e188}, and  \eqref{e189} in \eqref{e181}, we obtain
\begin{align}
\mathcal E_{7,6}(S_j(q))
&=\sum_{u,v\in\Z}(-1)^{u+v}
W_j\bigl(7(6u-1),7(6v-1)\bigr)q^{7(P_u+2P_v)}\notag\\
&=7^{2j}\sum_{u,v\in\Z}(-1)^{u+v}
W_j(6u-1,6v-1)q^{7(P_u+2P_v)}\notag\\
&\label{e190}=7^{2j}S_j(q^7).
\end{align}
where  $j=2,3$.

Employing \eqref{e173}, \eqref{e182}, and \eqref{e183} in \eqref{e190}, we obtain
\begin{align}\label{e191}
-7\mathcal E_{7,6}(H_2(q))&=-7\cdot7^4H_2(q^7)\\\intertext{and}
23\mathcal E_{7,6}(H_3(q))&\label{e192}=23\cdot7^6H_3(q^7).
\end{align}
Since \eqref{e191} and \eqref{e192} are exact identities in $\Z[[q]]$, we may cancel the nonzero integer factors $-7$ and $23$, respectively. Then, using \eqref{e173}, we obtain
$$
\mathcal E_{7,6}(H_j(q))=7^{2j}H_j(q^7)=7^{2j}f_7f_{14}M(q^7)^j,$$
where  $j=2,3$. This completes the proof of  \eqref{e174}.
\end{proof}

\section{Generating functions and divisibility properties}

In this section, we establish generating functions for $COS(n)$ and $CES(n)$ and prove their divisibility properties.

\begin{theorem}
We have
\begin{align}\label{e44}
\sum_{n\ge0}COS(n)q^n
&=\dfrac{L(q)-2L(q^2)}{f_1f_2}=\dfrac{2\phi(q)^4-\phi(-q)^4-1}{24f_1f_2}=\dfrac1{24}\left(
2\dfrac{f_2^{19}}{f_1^9f_4^8}
-\dfrac{f_1^7}{f_2^5}
-\dfrac1{f_1f_2}\right)\\\intertext{and}
\sum_{n\ge0}CES(n)q^n&\label{e45}=\dfrac{4L(q^2)}{f_1f_2}=\dfrac{1-E_2(q^2)}{6f_1f_2}.
\end{align}
\end{theorem}
\begin{proof}
From  \eqref{e2}, we obtain
\begin{align}
\sum_{n\ge0}COS(n)q^n
&=\left.
\dfrac{\partial}{\partial z}
\prod_{m\ge1}
\dfrac{1}
{\bigl(1-(zq)^{2m-1}\bigr)(1-q^{2m})^2}
\right|_{z=1}\notag\\
&=\dfrac1{f_1f_2}
\sum_{m\ge1}
\dfrac{(2m-1)q^{2m-1}}{1-q^{2m-1}}\notag\\
&\label{e46}=\dfrac{L(q)-2L(q^2)}{f_1f_2}, 
\end{align}
where the exponent of $z$ records the sum
of the odd parts.
Similarly, for two colors of each even part, we obtain
\begin{align}
\sum_{n\ge0}CES(n)q^n
&=\left.
\dfrac{\partial}{\partial z}
\prod_{m\ge1}
\dfrac{1}
{(1-q^{2m-1})\bigl(1-(zq)^{2m}\bigr)^2}
\right|_{z=1}\notag\\
&=\dfrac1{f_1f_2}
\sum_{m\ge1}\dfrac{4mq^{2m}}{1-q^{2m}}\notag\\
&\label{e47}=\dfrac{4L(q^2)}{f_1f_2}, 
\end{align}
where the exponent of $z$ records  the sum of the even parts.
From \eqref{e5}, \eqref{e9}, \eqref{e27}, and \eqref{e28}, we obtain
\begin{equation}\label{n4}
L(q)-2L(q^2)=\dfrac{M-1}{24}=\dfrac{2\phi(q)^4-\phi(-q)^4-1}{24}=\dfrac1{24}\left(
2\dfrac{f_2^{20}}{f_1^8f_4^8}
-\dfrac{f_1^8}{f_2^4}-1\right).
\end{equation}
Employing \eqref{n4} in \eqref{e46} and simplifying, we complete the proof of  \eqref{e44}.

Again, from \eqref{e21},  we note that 
\begin{equation}\label{n5}4L(q^2)=(1-E_2(q^2))/6.\end{equation} Employing \eqref{n5} in  \eqref{e47}, we complete the proof of \eqref{e45}.
\end{proof}

\begin{theorem}\label{e48}
For any  integer $n\ge 0$, we have
\begin{align}
\label{i1}&CES(n)=4\sum_{j=1}^{\lfloor n/2\rfloor}
\sigma(j)a(n-2j)\\\intertext{and}
\label{i2}&COS(n)=na(n)-4\sum_{j=1}^{\lfloor n/2\rfloor}
\sigma(j)a(n-2j).
\end{align}
Here $\lfloor x\rfloor$ denotes the greatest integer not exceeding $x$, and the sums are interpreted as $0$ when $n=0$ or $n=1$.
\end{theorem}
\begin{proof}
Employing \eqref{e2} and \eqref{e12} in \eqref{e47}, we obtain
\begin{equation}\label{e49}
\sum_{n\ge0}CES(n)q^n=4\left(\sum_{j\ge1}\sigma(j)q^{2j}\right)
\left(\sum_{m\ge0}a(m)q^m\right)=4\sum_{n\ge0}
\left(\sum_{j=1}^{\lfloor n/2\rfloor}
\sigma(j)a(n-2j)\right)q^n.
\end{equation}
Equating the coefficients of $q^n$ on both sides of  \eqref{e49}, we arrive at \eqref{i1}.  Employing  \eqref{e4} in \eqref{i1} and simplifying, we complete the proof of \eqref{i2}.
\end{proof}

The following corollary is an easy consequence of Theorem~\ref{e48}.
\begin{corollary}
For every nonnegative integer $n$, we have
\begin{align*}
&CES(n)\equiv0\pmod4\\\intertext{and}
&COS(4n)\equiv0\pmod4.
\end{align*}
\end{corollary}

\begin{theorem}
For any integer $n\ge 0$, we have
\begin{align}
	&\label{e62}CES(9n+1)\equiv0\pmod3,\\
	&\label{e63}CES(9n+8)\equiv0\pmod3, \\
	&\label{e60}COS(9n+8)\equiv0\pmod3,\\
	&\label{e59}COS(9n+7)\equiv0\pmod3,\\
&	\label{e58} COS(9n)\equiv0\pmod3,\\ \intertext{and}
&\label{e61}CES(9n)\equiv0\pmod3.
\end{align}
\end{theorem}
\begin{proof}
Employing  \eqref{e6}, \eqref{e14}, and \eqref{e51} in \eqref{e2}, we obtain
\begin{equation}\label{e64}
A(q)=\dfrac{\left(f_2^2/f_1\right)}{f_2^3}\equiv\dfrac{\psi(q)}{f_6}=\dfrac{f_9^2}{f_3f_{18}}+q\dfrac{\psi(q^9)}{f_6}\pmod3.
\end{equation}
Employing \eqref{e6},  \eqref{e16},  and \eqref{e14} in \eqref{e64}, we obtain
\begin{align}\label{e65}
&\mathcal E_{3,2}(A(q))\equiv0\pmod3\\\intertext{and}
&\label{e66}\mathcal E_{3,1}(A(q))\equiv\dfrac{\psi(q^3)}{f_2}=\dfrac{f_6^2}{f_3f_2}\equiv\dfrac{f_6}{f_3}f_2^2\pmod3.
\end{align}
Replacing $q$ by $q^2$ in \eqref{e50}, we obtain
\begin{equation}\label{e67}
L(q^2)\equiv L(q^6)+q^2f_6^8\pmod3.
\end{equation}
Employing \eqref{e64} and \eqref{e67} in \eqref{e45} and simplifying, we obtain
\begin{align}
\sum_{n\ge0}CES(n)q^n
\label{e68}&\equiv
\dfrac{f_9^2L(q^6)}{f_3f_{18}}
+q\dfrac{\psi(q^9)L(q^6)}{f_6}
+q^2\dfrac{f_6^8f_9^2}{f_3f_{18}}
+q^3\psi(q^9)f_6^7 \pmod 3.
\end{align}
Employing  \eqref{e6}, \eqref{e14}, \eqref{e16}, and \eqref{e13} in \eqref{e68}, we obtain
\begin{align}
\mathcal E_{3,1}\!\left(\sum_{n\ge0}CES(n)q^n\right)&\label{e69}\equiv\dfrac{\psi(q^3)L(q^2)}{f_2}\equiv\dfrac{f_6}{f_3}f_2^2L(q^2)\equiv-\dfrac{f_6}{f_3}D(f_2^2)\pmod3.\\
\intertext{Applying $\mathcal E_{3,2}$ to \eqref{e68}, only the third summand contributes. Therefore, using $f_2^3\equiv f_6\pmod3$, we obtain}
\mathcal E_{3,2}\!\left(\sum_{n\ge0}CES(n)q^n\right)&\label{e70}\equiv\dfrac{f_2^8f_3^2}{f_1f_6}=\dfrac{f_3^2}{f_1f_2}\dfrac{f_2^9}{f_6}\equiv\dfrac{f_3^2f_6^2}{f_1f_2}=f_3^2f_6^2A(q)\pmod3.
\end{align}
Employing \eqref{e16}, \eqref{e18}, and \eqref{e20} in \eqref{e69} and \eqref{e70}, and then simplifying by using \eqref{e65}, we obtain
\begin{align}\label{e71}
\sum_{n\ge0}CES(9n+1)q^n&\equiv-\dfrac{f_2}{f_1}\mathcal E_{3,0}\bigl(D(f_2^2)\bigr)\equiv0\pmod3,\\\intertext{and}
\sum_{n\ge0}CES(9n+8)q^n&\label{e72}\equiv f_1^2f_2^2\mathcal E_{3,2}(A(q))\equiv0\pmod3, 
\end{align}respectively. Now, \eqref{e62} and \eqref{e63} follow immediately from \eqref{e71} and \eqref{e72}, respectively. 
Replacing $n$ by $9n+8$ in  \eqref{e4} and simplifying using \eqref{e72} and \eqref{e65}, we arrive at  \eqref{e60}.

Employing \eqref{e18}, \eqref{e66}, and \eqref{e69} in \eqref{e4}, we obtain
\begin{align}
\mathcal E_{3,1}\!\left(\sum_{n\ge0}COS(n)q^n\right)
&=\mathcal E_{3,1}(DA(q))
-\mathcal E_{3,1}\!\left(\sum_{n\ge0}CES(n)q^n\right)\notag\\
&\equiv\dfrac{f_6}{f_3}f_2^2
+\dfrac{f_6}{f_3}D(f_2^2)\notag\\
\label{e73}&=\dfrac{f_6}{f_3}(1+D)f_2^2 \pmod3.
\end{align}
From \eqref{e17}, we note that
$$
\mathcal E_{3,2}\bigl((1+D)f_2^2\bigr)
=3(D+1)\mathcal E_{3,2}(f_2^2).
$$
Employing \eqref{e16} and \eqref{e20} in \eqref{e73}, we obtain
\begin{equation}\label{e74}
\sum_{n\ge0}COS(9n+7)q^n\equiv\dfrac{f_2}{f_1}\mathcal E_{3,2}\bigl((1+D)f_2^2\bigr)\equiv0\pmod3.
\end{equation} Now, \eqref{e59} follows easily from \eqref{e74}.

Again, from  \eqref{e2} and \eqref{e13}, we note that
\begin{equation}\label{e75}
\dfrac{DA(q)}{A(q)}=-\dfrac{Df_1}{f_1}-\dfrac{Df_2}{f_2}=L(q)+2L(q^2).
\end{equation}
Employing \eqref{e75} and \eqref{e9} in \eqref{e46}, we obtain
\begin{align}
\sum_{n\ge0}COS(n)q^n
&=A(q)\bigl(2L(q)-(L(q)+2L(q^2))\bigr)\notag\\
&=2A(q)L(q)-DA(q)\notag\\
&\label{e76}=2\dfrac{\phi(-q)L(q)}{f_1^3}-DA(q),
\end{align}
Using \eqref{e50} and \eqref{e14}, we obtain
\begin{equation}\label{e77}
\dfrac{L(q)}{f_1^3}\equiv\dfrac{L(q^3)+qf_3^8}{f_3}=\dfrac{L(q^3)}{f_3}+qf_3^7\pmod3.
\end{equation}
Also, from \eqref{e9}, we note that
\begin{equation}\label{e78}
\phi(-q)=\sum_{n\in\Z}(-1)^nq^{n^2}=\phi(-q^9)+\sum_{\substack{n\in\Z\\3\nmid n}}(-1)^nq^{n^2}.
\end{equation}
Employing \eqref{e13}, \eqref{e16}, \eqref{e18},  \eqref{e77}, and \eqref{e78} in \eqref{e76}, we obtain
\begin{align}
\mathcal E_{3,0}\!\left(\sum_{n\ge0}COS(n)q^n\right)
&\equiv2\mathcal E_{3,0}\!\left(\phi(-q^9)\dfrac{L(q^3)}{f_3}\right)\notag\\
&=2\phi(-q^3)\dfrac{L(q)}{f_1}\notag\\
&\label{e79}=2\phi(-q^3)D\!\left(\dfrac1{f_1}\right)\pmod3.
\end{align}
Employing \eqref{e16}, \eqref{e18}, and \eqref{e20} in \eqref{e79}, we obtain
\begin{equation}\label{e80}
\sum_{n\ge0}COS(9n)q^n\equiv2\phi(-q)\mathcal E_{3,0}\!\left(D\!\left(\dfrac{1}{f_1}\right)\right)\equiv0\pmod3.
\end{equation}Now, \eqref{e58} follows immediately from \eqref{e80}. Congruence  \eqref{e61} follows from \eqref{e58} and \eqref{e4}.
\end{proof}

\begin{theorem}
For any integer $n\ge 0$, we have
\begin{align}
&\label{e151}COS(25n+7)\equiv0\pmod5,\\
&\label{e152}COS(25n+22)\equiv0\pmod5,\\
&\label{e153}CES(25n+12)\equiv0\pmod5,\\\intertext{and}
&\label{e154}CES(25n+22)\equiv0\pmod5.
\end{align}
\end{theorem}

\begin{proof}
Employing \eqref{e2} and \eqref{e27} in \eqref{e44}, we obtain
\begin{equation}\label{e155}
\sum_{n\ge0}COS(n)q^n=\dfrac{MA(q)-A(q)}{24}\equiv A(q)-MA(q)\pmod5.
\end{equation}
Employing \eqref{e15}, \eqref{e112}, and \eqref{e124} in \eqref{e155}, we obtain
\begin{equation}\label{e156}
\sum_{n\ge0}COS(5n+2)q^n\equiv\mathcal E_{5,2}(A(q))-\mathcal E_{5,2}(MA(q))\equiv(3-M)H\pmod5, 
\end{equation}
where 
\begin{equation}\label{e150}
	H=f_1^3f_2^3.
\end{equation}
From \eqref{e18}, \eqref{e124}, and \eqref{e150}, we obtain
\begin{equation}\label{e157}
\mathcal E_{5,2}(DA(q))\equiv2\mathcal E_{5,2}(A(q))\equiv6H\pmod5.
\end{equation}
Employing \eqref{e156} and \eqref{e157} in \eqref{e4}, we obtain
\begin{align}
\sum_{n\ge0}CES(5n+2)q^n
&=\mathcal E_{5,2}(DA(q))-\sum_{n\ge0}COS(5n+2)q^n\notag\\
&\label{e158}\equiv(3+M)H\pmod5.
\end{align}

Next, for any  integers $r, s\ge0$, set
$$
T_r=\dfrac{r(r+1)}2
$$ and 
\begin{equation}\label{e159}
x=2r+1,\quad y=2s+1,\quad m=T_r+2T_s.
\end{equation}
Replacing $q$ by $q^2$ in \eqref{e11}, we obtain
\begin{equation}\label{e11e}
	f_2^3=\sum_{m\ge0}(-1)^m(2m+1)q^{m(m+1)}.
\end{equation}
Simplifying \eqref{e150} by using \eqref{e11}, \eqref{e11e}, and  \eqref{e159}, we obtain
\begin{equation}\label{e160}
H=\sum_{r,s\ge0}(-1)^{r+s}xyq^m.
\end{equation}
From \eqref{e13} and \eqref{e21}, we obtain
\begin{align}\label{e161}
&(8D+1)f_1^3=(1-24L(q))f_1^3=E_2(q)f_1^3\\\intertext{and}
&\label{e162}(4D+1)f_2^3=(1-24L(q^2))f_2^3=E_2(q^2)f_2^3.
\end{align}
Simplifying \eqref{e161} and \eqref{e162} by employing \eqref{e11} and \eqref{e159}, we obtain
\begin{align}\label{e163}
&E_2(q)f_1^3=\sum_{r\ge0}(-1)^r(2r+1)(8T_r+1)q^{T_r}=\sum_{r\ge0}(-1)^rx^3q^{T_r}\\\intertext{and}
&\label{e164}E_2(q^2)f_2^3=\sum_{s\ge0}(-1)^s(2s+1)(8T_s+1)q^{2T_s}=\sum_{s\ge0}(-1)^sy^3q^{2T_s},
\end{align}respectively.
Combining  \eqref{e11}, \eqref{e27},  \eqref{e163}  and \eqref{e164}, we obtain
\begin{equation}\label{e165}
MH=\bigl(2E_2(q^2)-E_2(q)\bigr)f_1^3f_2^3
=\sum_{r,s\ge0}(-1)^{r+s}xy(2y^2-x^2)q^m.
\end{equation}
Also, from \eqref{e159} note that
\begin{equation}\label{e166}
x^2+2y^2=(8T_r+1)+2(8T_s+1)=8m+3.
\end{equation}
If $xy\equiv 0\pmod 5$, then the corresponding summands in \eqref{e160} and \eqref{e165} vanish modulo $5$. Otherwise, each of $x^2$ and $y^2$ is congruent to $1$ or $4$ modulo $5$. Thus, by \eqref{e166}, we obtain the following table for $xy\not\equiv 0\pmod 5$, with all entries taken modulo $5$:
\begin{center}
\begin{tabular}{c c c c}
\toprule
$m$ & $(x^2,y^2)$ & $3+x^2-2y^2$ & $3(m-1)(m+1)^3$\\
\midrule
$0$ & $(1,1)$ & $2$ & $2$\\
$1$ & $(4,1)$ & $0$ & $0$\\
$2$ & $(1,4)$ & $1$ & $1$\\
$3$ & $(4,4)$ & $4$ & $4$\\
$4$ & impossible & -- & $0$\\
\bottomrule
\end{tabular}
\end{center}
Thus, from the above table, we always have
\begin{equation}\label{e167}
xy(3+x^2-2y^2)\equiv3(m-1)(m+1)^3xy\pmod5.
\end{equation}
Simplifying \eqref{e160} and \eqref{e165} and employing  \eqref{e167}, we obtain
\begin{align}
(3-M)H
&=\sum_{r,s\ge0}(-1)^{r+s}xy(3+x^2-2y^2)q^m\notag\\
&\equiv3\sum_{r,s\ge0}(-1)^{r+s}(m-1)(m+1)^3xyq^m\notag\\
&\label{e168}=3(D-1)(D+1)^3H \pmod5.
\end{align}
Employing \eqref{e168} in \eqref{e156}, we obtain
\begin{equation}\label{e169}
\sum_{n\ge0}COS(5n+2)q^n\equiv(3-M)H\equiv3(D-1)(D+1)^3H\pmod5.
\end{equation}
Similarly, from \eqref{e158} and \eqref{e168}, we obtain
\begin{equation}\label{e170}
\sum_{n\ge0}CES(5n+2)q^n\equiv(3+M)H\equiv3(2-(D-1)(D+1)^3)H\pmod5.
\end{equation}
Extracting the coefficients of $q^{5n+1}$ and $q^{5n+2}$ on both sides of \eqref{e169} and \eqref{e170}, respectively,  we obtain
\begin{align*}
COS(25n+7)&\equiv3(5n)(5n+2)^3[q^{5n+1}]H\equiv0\pmod5\\\intertext{and}
CES(25n+12)&\equiv3\bigl(2-(5n+1)(5n+3)^3\bigr)[q^{5n+2}]H
\equiv3(2-27)[q^{5n+2}]H\equiv0\pmod5.
\end{align*}
So the proofs of  \eqref{e151} and \eqref{e153} are complete.

Since $-2\equiv3\pmod5$ is a quadratic nonresidue modulo $5$, from \eqref{e166}, we obtain
\begin{equation}\label{e171}
m\equiv4\pmod5 \text{ if and only if }
\begin{cases}
x\equiv0\pmod5\\
y\equiv0\pmod5,
\end{cases}\text{if and only if }
\begin{cases}
r=5u+2\\
s=5v+2,
\end{cases}
\end{equation}
where $ u$ and $v$ are non-negative integers.

From \eqref{e159} and \eqref{e171}, we obtain
\begin{align}
\label{u1}T_{5u+2}&=\dfrac{(5u+2)(5u+3)}2=25T_u+3,\\
\label{u2}m&=25(T_u+2T_v)+9,\\
\label{u3}(-1)^{r+s}&=(-1)^{5u+5v+4}=(-1)^{u+v},\\\intertext{and}
\label{u4}xy&=(10u+5)(10v+5)=25(2u+1)(2v+1). 
\end{align}
Substituting \eqref{u1}, \eqref{u2}, \eqref{u3}, and \eqref{u4} in \eqref{e160} and using \eqref{e11}, we obtain
\begin{align}
\mathcal E_{5,4}(H)
&=25\sum_{u,v\ge0}(-1)^{u+v}(2u+1)(2v+1)
q^{(25(T_u+2T_v)+9-4)/5}\notag\\
&=25q\left(\sum_{u\ge0}(-1)^u(2u+1)q^{5T_u}\right)
\left(\sum_{v\ge0}(-1)^v(2v+1)q^{10T_v}\right)\notag\\
&\label{e172}=25qf_5^3f_{10}^3.
\end{align}
Now, by \eqref{e172} it is easily seen that
$$
[q^{5n+4}]H=25[q^n]\bigl(qf_5^3f_{10}^3\bigr)\equiv0\pmod5.
$$
Extracting the coefficients of $q^{5n+4}$ on both sides of \eqref{e169} and \eqref{e170}, we obtain
\begin{align*}
COS(25n+22)
&\equiv3(5n+3)(5n+5)^3[q^{5n+4}]H\equiv0\pmod5\\\intertext{and}
CES(25n+22)
&\equiv3\bigl(2-(5n+3)(5n+5)^3\bigr)[q^{5n+4}]H
\equiv6[q^{5n+4}]H\equiv0\pmod5,
\end{align*}
respectively. Hence, the proofs of  \eqref{e152} and \eqref{e154} are complete.
\end{proof}
\begin{theorem}
For every nonnegative integer $n$, we have
\begin{align}\label{e194}
COS(49n+43)&\equiv0\pmod7\\\intertext{and}
CES(49n+43)&\label{e195}\equiv0\pmod7.
\end{align}
\end{theorem}

\begin{proof}
Employing \eqref{e2} and \eqref{e27} in \eqref{e44}, we obtain
\begin{equation}\label{e196}
\sum_{n\ge0}COS(n)q^n
=\dfrac{MA(q)-A(q)}{24}
\equiv5(MA(q)-A(q))\pmod7.
\end{equation}
Employing \eqref{e15}, \eqref{e113}, and \eqref{e125} in \eqref{e196}, we obtain
\begin{align}
\sum_{n\ge0}COS(7n+1)q^n
&\equiv5\bigl(\mathcal E_{7,1}(MA(q))-\mathcal E_{7,1}(A(q))\bigr)\notag\\
&\equiv5\bigl(4f_1f_2M^3-f_1f_2M^2\bigr)\notag\\
&=f_1f_2M^2(20M-5)\notag\\
&\label{e197}\equiv f_1f_2M^2(2+6M)\pmod7.
\end{align}
Also, from \eqref{e18} and \eqref{e125}, we obtain
\begin{equation}\label{e198}
\mathcal E_{7,1}(DA(q))\equiv\mathcal E_{7,1}(A(q))\equiv f_1f_2M^2\pmod7.
\end{equation}
Employing \eqref{e197} and \eqref{e198} in \eqref{e4}, we obtain
\begin{align}
\sum_{n\ge0}CES(7n+1)q^n
&=\mathcal E_{7,1}(DA(q))-\sum_{n\ge0}COS(7n+1)q^n\notag\\
&\equiv f_1f_2M^2-f_1f_2M^2(2+6M)\notag\\
&=f_1f_2M^2(-1-6M)\notag\\
&\label{e193}\equiv f_1f_2M^2(6+M)\pmod7.
\end{align}
Again, from \eqref{e173} and \eqref{e174}, we obtain
\begin{align}\label{e199}
\mathcal E_{7,6}(f_1f_2M^2)
&=7^4f_7f_{14}M(q^7)^2\equiv0\pmod7\\\intertext{and}
\mathcal E_{7,6}(f_1f_2M^3)
&\label{e200}=7^6f_7f_{14}M(q^7)^3\equiv0\pmod7.
\end{align}
Employing \eqref{e20}, \eqref{e199}, and \eqref{e200} in \eqref{e197} and \eqref{e193}, we obtain
\begin{align}
\label{l1}\sum_{n\ge0}COS(49n+43)q^n
&\equiv2\mathcal E_{7,6}(f_1f_2M^2)
+6\mathcal E_{7,6}(f_1f_2M^3)\equiv0\pmod7\\\intertext{and}
\label{l2}\sum_{n\ge0}CES(49n+43)q^n
&\equiv6\mathcal E_{7,6}(f_1f_2M^2)
+\mathcal E_{7,6}(f_1f_2M^3)\equiv0\pmod7,
\end{align}
Now \eqref{e194} and \eqref{e195} follow easily from \eqref{l1} and \eqref{l2}, respectively. 
\end{proof}
\section{Concluding remarks}

Having established the above congruences modulo $3$, $5$, and $7$, it is natural to ask whether corresponding congruences hold modulo higher powers. Based on numerical evidence, we make the following conjecture.

\medskip
\noindent\textbf{Conjecture.}
Let $p\in\{3,5,7\}$ and let $\alpha$ be a positive integer. Then, for every nonnegative integer $n$,
\begin{equation}
COS\left(p^{2\alpha}n+\dfrac{7p^{2\alpha}+1}{8}\right)
\equiv
CES\left(p^{2\alpha}n+\dfrac{7p^{2\alpha}+1}{8}\right)
\equiv0\pmod{p^\alpha}.
\end{equation}

Since $p$ is odd, $p^{2\alpha}\equiv1\pmod8$, so $\dfrac{7p^{2\alpha}+1}{8}$ is an integer. We note that the case $\alpha=1$ gives precisely the congruences
proved above for $p=3$, $5$, and $7$.

\section*{\bf Declarations}
\noindent\textbf{Funding}: This research did not receive funding.\\
\noindent{\bf Author Contributions.} Both authors contributed equally to this work.

\noindent{\bf Conflict of Interest.} The authors declare that there is no conflict of interest regarding the publication of this paper.

\noindent{\bf Human and Animal Rights.} The authors declare that there is no research involving human participants or animals in the context of this paper.	

\noindent{\bf Data Availability Statement.} Data sharing is not applicable to this paper as no datasets were generated or analyzed during the current study.


\begin{thebibliography}{99}

\bibitem{AndrewsDastidar2026}
Andrews, G. E. and Ghosh Dastidar, M.:
$p(5n+4)$ again.
\emph{Ramanujan J.} \textbf{69} (2026), Article 26.

\bibitem{BS}
Bardhan, G. and Saikia, N.: Proofs of the conjectures on
$SOME(n)$ and $DSOME(n)$ functions related to integer
partitions. Preprint, arXiv:2608.08524 (2026).

\bibitem{BaruahGogoi2026}
Baruah, N. D. and Gogoi, P.:
Difference of the sum of all odd and even overlined parts of overpartitions.
Preprint, arXiv:2609.21361 (2026).

\bibitem{Berndt2006}
Berndt, B. C.:
\emph{Number Theory in the Spirit of Ramanujan}.
American Mathematical Society, Providence, RI (2006).

\bibitem{Chan2010}
Chan, H.-C.:
Ramanujan's cubic continued fraction and an analog of his ``most beautiful identity''.
\emph{Int. J. Number Theory} \textbf{6} (2010), 673--680.

\bibitem{ChenLin2009}
Chen, W. Y. C. and Lin, B. L. S.:
Congruences for the number of cubic partitions derived from modular forms.
Preprint, arXiv:0910.1263 (2009).

\bibitem{euler1748introductio}
Euler, L.:
\emph{Introductio in Analysin Infinitorum}, Vol.~1.
M.-M. Bousquet, Lausanne (1748).

\bibitem{GarvanSarma2026}
Garvan, F. G. and Sarma, R.:
Sum of parts in overpartitions and partitions without repeated odd parts.
Preprint, arXiv:2606.14680 (2026).

\bibitem{Jacobi1829}
Jacobi, C. G. J.:
\emph{Fundamenta Nova Theoriae Functionum Ellipticarum}.
Fratres Borntraeger, Regiomonti (1829).

\bibitem{RSE}
Ramanujan, S.: On certain arithmetical functions.
\emph{Trans. Cambridge Philos. Soc.} \textbf{22} (1916), 159--184.

\bibitem{Ramanujan}
Ramanujan, S.:
Congruence properties of partitions.
\emph{Math. Z.} \textbf{9} (1921), 147--153.



\bibitem{tifr}
Ramanujan, S.: \emph{Notebooks}, 2 vols.
Tata Institute of Fundamental Research, Bombay (1957).

\bibitem{lost}
Ramanujan, S.:
\emph{The Lost Notebook and Other Unpublished Papers}.
Narosa Publishing House, New Delhi (1988).

\bibitem{Xiong2011}
Xiong, X.:
The number of cubic partitions modulo powers of $5$ (in Chinese).
\emph{Sci. Sin. Math.} \textbf{41} (2011), 1--15;
English preprint: arXiv:1004.4737.

\end{thebibliography}
\end{document}